\documentclass[5p,twocolumn,times]{elsarticle}

\usepackage{algorithm}
\usepackage{algorithmic}
\usepackage{amsmath}
\usepackage{amssymb}
\usepackage{color}
\usepackage{graphicx}
\usepackage{listings}
\usepackage{subfig}
\usepackage{hyperref}
\usepackage{xspace}
\usepackage{textcomp} 

\definecolor{blue}{rgb}{0.098,0.357,0.675}
\definecolor{green}{rgb}{0.5,0.75,0.0}

\lstdefinelanguage{CUDA}[]{C++}{
    morekeywords={__global__,__device__,__shared__,__syncthreads,threadIdx,blockIdx,float3,float4,rsqrtf},
}
\newcommand{\dns}{\textsc{dns}\xspace}
\newcommand{\cfd}{\textsc{cfd}\xspace}
\newcommand{\fmm}{\textsc{fmm}\xspace}
\newcommand{\fft}{\textsc{fft}\xspace}
\newcommand{\mac}{\textsc{mac}\xspace}
\newcommand{\mpi}{\textsc{mpi}\xspace}
\newcommand{\orb}{\textsc{orb}\xspace}
\newcommand{\LET}{\textsc{let}\xspace}
\newcommand{\cpu}{\textsc{cpu}}
\newcommand{\gpu}{\textsc{gpu}}

\newcommand{\nvidia}{N\textsc{vidia}\xspace}
\newcommand{\tsubame}{\textsc{tsubame}{\footnotesize -2.0}\xspace}
\newcommand{\ccby}{\textsc{cc-by}{\footnotesize -3.0}\xspace}
\newcommand{\exafmm}{\texttt{exaFMM}\xspace}
\newcommand{\PM}{\textsc{p}\texttwooldstyle\textsc{m}\xspace} 
\newcommand{\LP}{\textsc{l}\texttwooldstyle\textsc{p}\xspace} 

\newcommand{\MM}{\textsc{m}\texttwooldstyle\textsc{m}\xspace} 
\newcommand{\ML}{\textsc{m}\texttwooldstyle\textsc{l}\xspace} 
\newcommand{\LL}{\textsc{l}\texttwooldstyle\textsc{l}\xspace}  
\newcommand{\PP}{\textsc{p}\texttwooldstyle\textsc{p}\xspace} 

\begin{document}

\begin{frontmatter}

\title{Petascale turbulence simulation using a highly parallel fast multipole method on GPUs}

\author[bu]{Rio Yokota}
\ead{yokota@bu.edu}

\author[bu]{L.~A.~Barba\corref{lab}}
\ead{labarba@bu.edu}

\author[uec]{Tetsu Narumi}
\ead{narumi@cs.uec.ac.jp}

\author[keio]{Kenji Yasuoka}
\ead{yasuoka@mech.keio.ac.jp}

\address[bu]{Dept.\ of Mechanical Engineering, Boston University, Boston MA 02215}
\address[uec]{Dept.\ of Computer Science, University of Electro-Communications, Tokyo, Japan}
\address[keio]{Dept.\ of Mechanical Engineering, Keio University, Yokohama, Japan}
\cortext[lab]{Correspondence to:  110 Cummington St, Boston MA 02215, (617) 353-3883,
\href{mailto:labarba@bu.edu}{labarba@bu.edu}}

\begin{abstract}
This paper reports large-scale direct numerical simulations of homogeneous-isotropic fluid turbulence, achieving sustained performance of 1.08 petaflop/s on \gpu\ hardware using single precision. The simulations use a vortex particle method to solve the Navier-Stokes equations, with a highly parallel fast multipole method (\fmm) as numerical engine, and match the current record in mesh size for this application, a cube of  $4096^3$ computational points solved with a spectral method.
The  standard numerical approach used in this field is the pseudo-spectral method, relying on the \fft algorithm as numerical engine. The particle-based simulations presented in this paper  quantitatively match the kinetic energy spectrum obtained with a pseudo-spectral method, using a trusted code. In terms of parallel performance, weak scaling results show the \fmm-based vortex method achieving 74\% parallel efficiency on 4096 processes (one \gpu\ per \mpi process, 3 \gpu s per node of the \tsubame system). The \fft-based spectral method is able to achieve just 14\% parallel efficiency on the same number of \mpi processes (using only \cpu\ cores), due to the all-to-all communication pattern of the \fft algorithm.
The calculation time for one time step was 108 seconds for the vortex method and 154 seconds for the spectral method, under these conditions. Computing with 69 billion particles, this work exceeds by an order of magnitude the largest vortex-method calculations to date.
\end{abstract}

\end{frontmatter}

\section{Introduction}

In the history of using computer simulation as a research tool to study the physics of turbulence, the dominant approach has been to use spectral methods. Direct numerical simulation (\dns) was introduced  as a means to check the validity of turbulence theories \emph{directly} from the equations of fluid dynamics \cite[p.~361]{Davidson2011}. The idea that important features of turbulence are \emph{universal} encouraged researchers to study the simplest of geometries, a periodic cubic volume of homogeneous, isotropic turbulent fluid. In this case, the simplicity and efficiency of a Fourier-spectral method cannot be matched. 
The largest direct numerical simulation of isotropic turbulence was conducted by Ishihara et al.~\cite{IshiharaETal2007,YokokawaETal2002} using $4096^3$ grid points at a maximum Taylor-microscale Reynolds number $R_\lambda =1200$. This record-breaking computation was done on \textsl{Earth Simulator}---a large vector machine with crossbar switch interconnect that can efficiently perform large-scale \fft. The successive generations of supercomputers have not been so \fft-friendly, and this record has not been surpassed even though the peak performance of supercomputers has increased nearly 50-fold since then. The record was matched for the first time in the U.S. by Donzis et al.~\cite{DonzisETal2008}, running on 16 thousand \cpu\ cores of the Ranger supercomputer in Texas, a Linux-cluster supercomputer with 16-core nodes.

Future high-performance computing systems will have ever more nodes, and ever more cores per node, but will probably not be equipped with the bandwidth required by many popular algorithms to transfer the necessary data at the optimal rate. This situation is detrimental to parallel scalability. Therefore, it is becoming increasingly important to consider alternative algorithms that may achieve better sustained performance on these extremely parallel machines of the future.

In most standard methods of incompressible \cfd, the greatest fraction of the calculation runtime is spent solving a Poisson equation. Equations of this type can be efficiently solved by means of an \fft-based algorithm, a sparse linear solver, or a fast multipole method (\fmm) \cite{GreengardRokhlin1987}. For the sake of our argument, we will not differentiate between \fft-based Poisson solvers and pseudo-spectral methods because they both rely on \fft. The fast multipole method has not gained popularity due to the fact that it is substantially slower---depending on implementations, at least an order of magnitude slower---than \fft and multigrid solvers when compared using a small \cpu\ cluster. We aim to show with our ongoing research that the relative performance of \fmm improves as one scales to large numbers of processes using \gpu\ acceleration.

The highly scalable nature of the \fmm algorithm, among other features, makes it a top contender in the algorithmic toolbox for exascale systems. One point of evidence for this argument is given by the Gordon Bell Prize winner of 2010, which achieved 0.7 petaflop/s with an \fmm algorithm on 200k cores of the Jaguar supercomputer at Oak Ridge National Laboratory \cite{RahimianETal2010}. In the previous year, \fmm also figured prominently at the Supercomputing Conference, with a paper among the finalists for the Best Paper award \cite{LashukETal2009} and the Gordon Bell prize in the price/performance category going to work with hierarchical $N$-body methods on \gpu\ architecture \cite{HamadaETal2009a}.

The \fmm algorithm is well adapted to the architectural features of \gpu s, which is an important consideration given that \gpu s are likely to be a dominant player as we move towards exascale. The work presented in 2009---winner in the price/performance category in great measure thanks to the ingenious and resourceful system design using gaming hardware---reported (at the time of the conference) 80 teraflop/s \cite{HamadaETal2009a}. That work progressed to an Honorable Mention in the 2010 list of awardees, with a performance of 104 teraflop/s \cite{HamadaNitadori2010}.

At the level of the present work, where we present a 1.08 petaflop/s (single precision) calculation of homogeneous isotropic turbulence, the \fmm moves firmly into the arena of \emph{petascale} \gpu\ computing. The significance of this advance is that we are now in a range where the \fmm algorithm shows its true capability. The excellent scalability of \fmm using over 4000 \gpu s is an advantage over the dominant \fft-based algorithms. Showcasing the \fmm in the simulation of homogeneous isotropic turbulence is especially fitting, given that a years-old record there remains unchallenged. We not only match the grid size of the world record, but also demonstrate that using the \fmm as the underlying algorithm enables isotropic turbulence simulations that scale to thousands of \gpu s. Given the severe bottleneck imposed by the all-to-all communication pattern of the \fft algorithm, this is not possible with pseudo-spectral methods in current hardware.

\section{Integration of \fmm into turbulence simulations}
The \fmm is used here as the numerical engine in a vortex particle method, which is an approach to solve the Navier-Stokes equations using a the vorticity formulation of momentum conservation and a particle discretization of vorticity. This is not a standard approach for the simulation of turbulence and the vortex method is not yet trusted for this application. For this reason, we have made efforts to document the validation of our vortex-method code and, in a separate publication \cite{YokotaBarba2012b}, compared with a trusted spectral-method code. We looked at various turbulence statistics, including higher-order velocity derivatives, and performed a parameter study for the relevant numerical parameters of the vortex method. That work provides evidence that the vortex method is an adequate tool for direct numerical simulation of fluid turbulence, while in the present work we focus on the performance aspects. For completeness, this sections gives a brief overview of the numerical methods.

\subsection{Vortex method}

The vortex method \cite{CottetKoumoutsakos2000} is a particle-based approach for fluid dynamics simulations. The particle discretization results in the continuum physics being solved as an $N$-body problem. Therefore, the hierarchical $N$-body methods that extract the full potential of \gpu s can be used for the simulation of turbulence. Unlike other particle-based solvers for fluid dynamics, e.g., smoothed particle hydrodynamics \cite{GingoldMonaghan1977}, the vortex method is especially well suited for computing turbulent flows, because the vortex interactions seen in turbulent flows are precisely what it calculates with the vortex particles. 

In the vortex method, the Navier-Stokes equation is solved in the velocity-vorticity formulation, using a particle discretization of the vorticity field. The velocity is calculated using the following equation, representing the Biot-Savart law of fluid dynamics:
\begin{equation}
\mathbf{u}_i=\sum_{j=1}^{N}\boldsymbol{\gamma}_j\times\nabla Gg_\sigma.
\label{eq:biotsavart}
\end{equation}
Here, $\boldsymbol{\gamma}_{j}$ is the strength of vortex particles, $G=1/4\pi r_{ij}$ is the Green's function for the Laplace equation and
\begin{equation}
g_\sigma=\mathrm{erf}\left(\sqrt{\frac{r_{ij}^{2}}{2\sigma_{j}^{2}}}\right)-\sqrt{\frac{4}{\pi}}\sqrt{\frac{r_{ij}^{2}}{2\sigma_{j}^{2}}}\exp\left(-\frac{r_{ij}^{2}}{2\sigma_{j}^{2}}\right)
\label{eq:cutoff}
\end{equation}
is the cutoff function, with $r$ the distance between the interacting particles, and $\sigma$ the standard deviation of the Gaussian function. The Navier-Stokes system is solved by a simultaneous update of the particle positions to account for convection, of the particle strengths to account for vortex stretching, and of the particle width to account for diffusion. The equation used to calculate the stretching term, $\boldsymbol{\omega}\cdot\nabla\mathbf{u}$, is:
\begin{equation}
\frac{{\rm d}\boldsymbol{\gamma}_i}{{\rm d} t}=\sum_{j=1}^{n}\nabla(\boldsymbol{\gamma}_j\times\nabla Gg_\sigma)\cdot\boldsymbol{\gamma}_i,
\label{eq:stretching}
\end{equation}
which was obtained by substituting the Biot-Savart equation (\ref{eq:biotsavart}) for $\mathbf{u}$, and using the discrete form of the vorticity. Finally, the diffusion update is calculated according to
\begin{equation}
\frac{{\rm d} \sigma^2}{{\rm d} t}=2\nu.
\label{eq:csm}
\end{equation}
We perform a radial basis function interpolation for reinitialized Gaussian distributions to ensure the convergence of the diffusion calculation \cite{Barba2004}.
Equations (\ref{eq:biotsavart}) and (\ref{eq:stretching}) are $N$-body interactions, and are evaluated using the \fmm. We use a highly parallel \fmm library for \gpu s developed in our group, called \texttt{exaFMM}, which is available under an open-source license and is described further below.

\subsection{Spectral method}\label{s:spectral}

\begin{figure}
\centering
\includegraphics[width=1.0\linewidth]{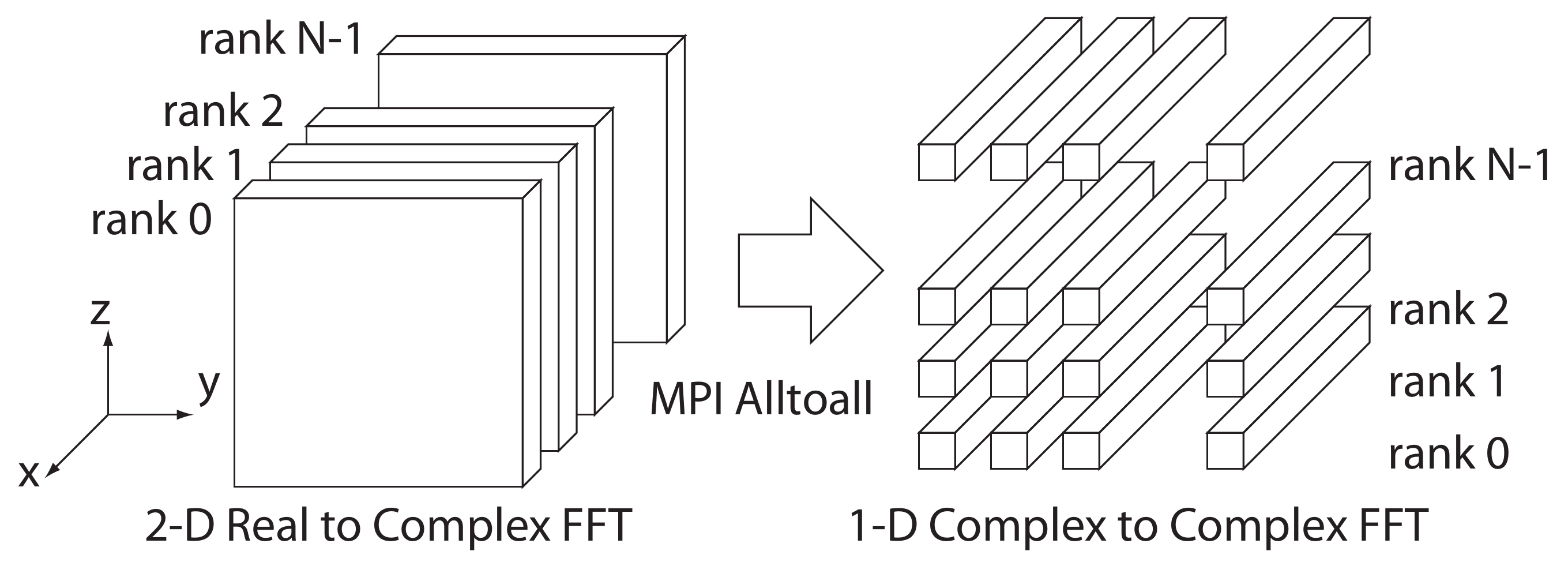}
\caption{Two stage parallel \fft in spectral method calculation}
\label{fig:spectral_fft}
\end{figure}

For the purpose of comparing with a spectral method, we used a code for homogeneous isotropic turbulence developed and used at the Center for Turbulence Research of Stanford University  \cite{ChumakovETal2009}. The code is called \texttt{hit3d} and is available freely for download.\footnote{In Google Code at \url{http://code.google.com/p/hit3d/}} It uses a spectral Galerkin method in primitive-variable formulation with pseudo-spectral methods to compute the convolution sums. For the \fft, it relies on the \texttt{fftw} library and it provides parallel capability for \cpu\ clusters using \mpi.  Parallelization of the \fft is accomplished by a domain-decomposition approach, illustrated by the schematic of Figure \ref{fig:spectral_fft}. Domain decomposition is applied in two spatial directions, resulting in first ``slabs'' then ``pencil'' sub-domains that are assigned to each \mpi process.

The initial conditions for our runs were generated by \texttt{hit3d}. The vortex method used the same initial condition by first calculating the vorticity field in physical space, and then using radial basis function interpolation to obtain the vortex strengths. This is different from our other publication, where we studied the accuracy of turbulence simulations using the \fmm-based vortex method on \gpu s, looking at high-order turbulence statistics \cite{YokotaBarba2012b}. For that case, the initial condition provided by \texttt{hit3d}, which has a fully developed energy spectrum, was not suitable for our validation exercise looking at the time evolution of the velocity derivative skewness and flatness. For this reason, we constructed initial conditions in Fourier space as a solenoidal isotropic velocity field with random phases and a prescribed energy spectrum. This initial velocity field had a Gaussian distribution and satisfied the incompressibility condition.

\section{Parallel fast multipole method on GPUs} 
It is common that algorithms with low complexity (sparse linear algebra, \fft) have low arithmetic intensity, while algorithms with high arithmetic intensity (dense linear algebra) tend to have high complexity. The \fmm possesses a rare combination of $\mathcal{O}(N)$ complexity and an arithmetic intensity that is even higher than \textsc{dgemm} \cite{YokotaBarba2012a}. Although this may seem like a great combination, it also implies that there is a large constant in front of the $\mathcal{O}(N)$ scaling, which results in a larger time-to-solution compared to other $\mathcal{O}(N)$ or $\mathcal{O}(N\log N)$ methods like multigrid methods and \fft. However, as arithmetic operations become cheaper compared to data movement in terms of both cost and energy, the large asymptotic constant of the \emph{arithmetic} complexity becomes less of a concern.

\begin{figure*}[t]
\centering
\includegraphics[width=0.95\linewidth]{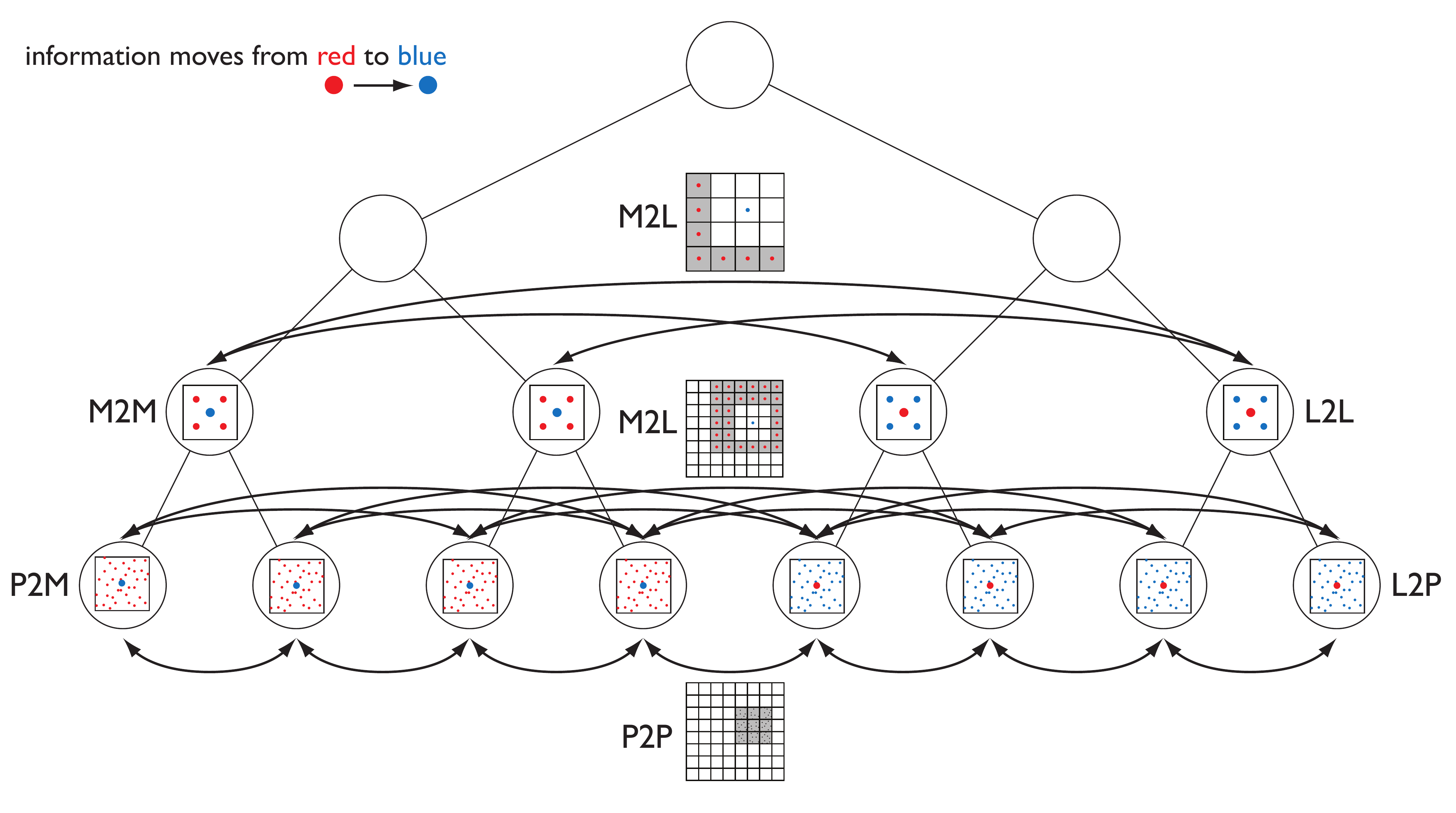}
\caption{Illustration of the flow of an \fmm calculation, with the naming convention used for the various computational kernels:  particle-to-multipole or \PM; multipole-to-multipole or \MM; multipole-to-local or \ML; local-to-local or \LL; local-to-particle or \LP; particle-to-particle or \PP (all explained in the text).}
\label{fig:kernels}
\end{figure*}

\subsection{Brief overview of the fast multipole method}

The calculation of velocity by means of the Biot-Savart equation (\ref{eq:biotsavart}) and the computation of the stretching term (\ref{eq:stretching})  both result in an $N$-body problem, which has a complexity of $\mathcal{O}(N^2)$ for $N$ particles, if solved directly. The fast multipole method \cite{GreengardRokhlin1987} reduces the complexity to $\mathcal{O}(N)$ by clustering source and target particles, and using series expansions that are valid far (multipole expansion) or near (local expansion) a point. The far/near relationships between points in the domain and the interactions between clusters are determined by means of a tree structure. The domain is hierarchically divided, and different sub-domains are associated with branches in the tree, a graphical representation of which is shown in Figure \ref{fig:kernels}. The algorithm proceeds as follows. First, the strengths of the particles (e.g., charges or masses) are transformed to multipole expansions at the leaf cells of the tree (known as the particle-to-multipole, or \PM kernel). Then, the multipole expansions of the smaller cells are translated to the center of larger cells in the tree hierarchy recursively and added (multipole-to-multipole, or \MM kernel). Subsequently, the multipole expansions are transformed into local expansions for all pairs of well-separated cells (multipole-to-local, or \ML kernel), and then to local expansions at the center of smaller cells recursively (local-to-local, or \LL kernel). Finally, the local expansions at leaf cells are used to compute the effect of the far field on each target particle. Since \ML operations can only be performed for well-separated cells, the direct neighbors at the finest level of the tree interact directly via the original equation (particle-to-particle, or \PP kernel). The current implementation of the \texttt{exaFMM} code uses expansions in spherical harmonics of the Laplace Green's function \cite{ChengETal1999}. Details of the extension of Laplace kernels to Biot-Savart and stretching kernels and of the implementation on \gpu s can be found in previous publications \cite{YokotaETal2009,YokotaBarba2010}.

\subsection{Tree partitioning} \label{sse:orb}

When parallelizing hierarchical $N$-body algorithms, the fact that the particle distribution is dynamic makes it impossible to precompute the decomposition and to balance the work-load or communication \textit{a priori}. Warren and Salmon \cite{WarrenSalmon1992} developed a parallel algorithm for decomposing the domain into recursive subdomains using the method of orthogonal recursive bisection (\orb). The resulting process can be thought of as a binary tree, which splits the domain into subdomains with equal number of particles at every bisection.
Another popular technique for partitioning tree structures is to use Morton ordering \cite{WarrenSalmon1993}, where bits representing the particle coordinates are interleaved to form a unique key that maps to each cell in the tree. Following the Morton index monotonically will take the form of a space-filling curve in the shape of a ``Z".  Partitioning the list of Morton indices equally assures that each partition will contain an equal number of cells, regardless of the particle distribution.

\begin{figure*}[t]
\centering
\includegraphics[width=0.95\linewidth]{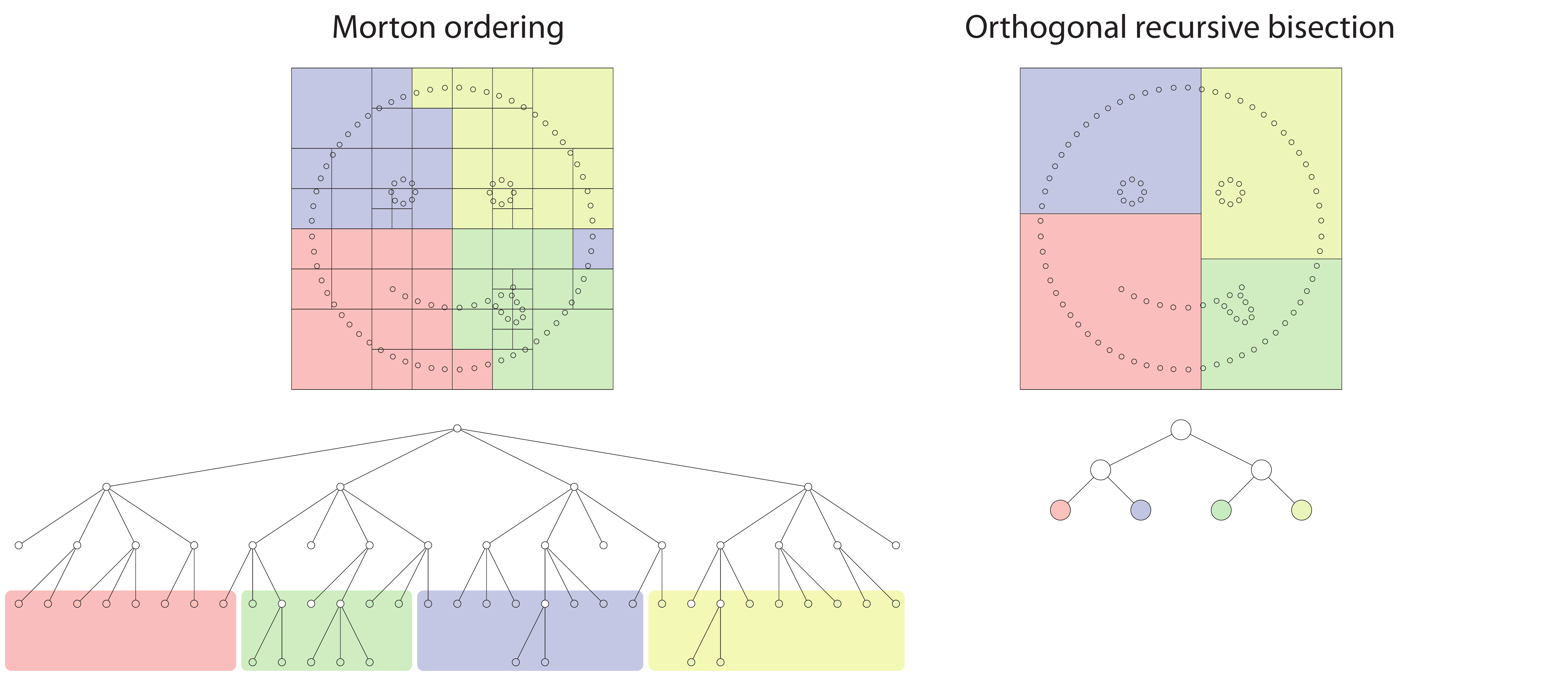}
\caption{Example of partitioning with Morton ordering (left) or with orthogonal recursive bisection, \orb (right) for a non-uniform particle distribution. The \orb partition always gives a balanced tree.}
\label{fig:partitioning}
\end{figure*}

For an adaptive tree, a na\"{\i}ve implementation of the Morton ordering could result in large communication costs if the ``Z" is split in the wrong place, as illustrated in Figure \ref{fig:partitioning}. Sundar et al.\ \cite{SundarETal2008} proposed a bottom-up coarsening strategy that ensures a clean partition at the coarse level while using Morton ordering. On the other hand, an \orb always partitions the domain into a well balanced binary tree, since the number of particles is always equal on each side of the bisection. Therefore, \orb is advantageous from the point of view of load-balancing.

The main difference between using Morton ordering and using an \orb is the shape of the resulting subdomains and the corresponding tree structure that connects them, as shown in Figure \ref{fig:partitioning}. Morton ordering is done on cubic octrees and the corresponding tree structure can become highly non-uniform depending on the particle distribution. Conversely, an \orb always creates a balanced tree structure but the sub-domain shapes are rectangular and non-uniform.

A difficulty when applying \orb within the framework of conventional \fmm is that the construction of cell-cell interaction lists depends on the cells being cubic (and not rectangular). The dual tree traversal described in the following subsection allows the use of rectangular cells, while enabling cell-cell interactions in the \fmm with minimum complications. By using \orb along with the dual tree traversal, we obtain a partitioning+traversal mechanism that allows perfect load balancing for highly non-uniform particle distributions, while retaining the $\mathcal{O}(N)$ complexity of the \fmm. Another advantage of this technique is that global indexing of cells is no longer necessary. Having a global index becomes an issue when solving problems of large size. The maximum depth of octrees that a 64-bit integer can handle is 21 levels ($2^{63}=8^{21}$). For highly adaptive trees with small number of particles per cell, this limit will be reached quite easily, resulting in integer-overflow with a global index. Circumventing this problem by using multiple integers for index storage will require much more work when sorting, so this is not a desirable solution. The dual tree traversal allows the entire \fmm to be performed without the use of global indices, and is an effective solution to this problem.

Our current partitioning scheme is an extension of the \orb, which allows multi-sections instead of bisections. Bisecting the domain involves the calculation of the median of the particle distribution for a given direction, and doing this recursively in orthogonal directions $(x,y,z,x,y,\ldots)$ is what constitutes the ``orthogonal recursive bisection". Therefore the extension from bisection to multi-section can be achieved by simply providing a mechanism to search for something other than the median. We developed a parallel version of the ``$n^{\text{th}}$-element" algorithm. Finding the ``$n^{\text{th}}$-element" is much faster than any sorting algorithm, so our technique is much faster than any method that requires sorting. Always searching for the $N/2$-th element will reduce to the original algorithm based on bisections, but searching for  the $3N/7$-th element will enable the domain to be split between 3 and 4 processes, for example. Therefore, our recursive multisection allows efficient partitioning when the number of processes is not a power of two.

\begin{figure}[t]
\centering
\includegraphics[width=1.0\linewidth]{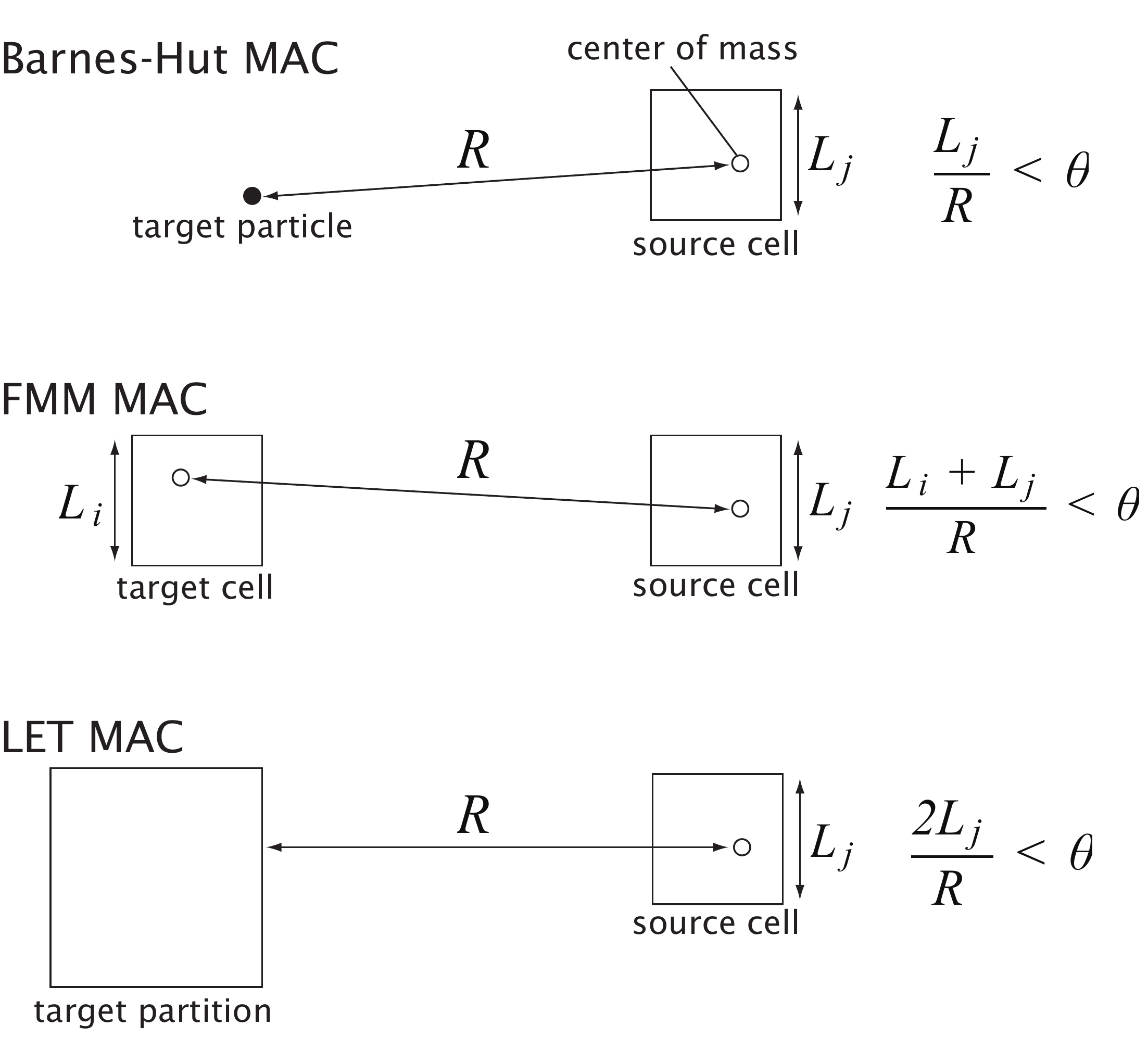}
\caption{Definition of three types of multipole acceptance criterion (\mac). The Barnes-Hut \mac is the simplest version from the original treecode. The \fmm \mac is defined between two cells rather than between a cell and particle. The \LET \mac is used for the construction of the local essential tree (\LET), and is defined between a cell and a partition on a remote process.}
\label{fig:mac}
\end{figure}

\subsection{Dual tree traversal} \label{sse:dual}

The dual tree traversal  enables cell-cell interactions in the $\mathcal{O}(N\log N)$ Barnes-Hut treecode \cite{BarnesHut1986} framework, thus turning it into a $\mathcal{O}(N)$ algorithm \cite{Dehnen2002}. We give a detailed explanation of the dual tree traversal in a previous publication that focused on hybrid treecode \& \fmm techniques \cite{YokotaBarba2012a}. It has several advantages compared to explicit construction of cell-cell interaction lists, the more common approach in the \fmm framework. Firstly, it is the simplest and most general way to perform the task of finding all combinations of cell-cell interactions. It is simple in the sense that it does not contain multiple loops for finding ``the parent's neighbor's children that are non-neighbors of the current cell," as is the case with the explicit construction. It is general in the sense that the property ``well-separated'' can be defined flexibly, instead of the rigid definition of ``non-neighboring cells" used traditionally in the \fmm. It is common in treecodes to define the well-separated property more precisely and adaptively by introducing the concept of multipole acceptance criterion (\mac) \cite{SalmonWarren1994}.

We give a graphical description of the different types of \mac used in our method in Figure \ref{fig:mac}. The simplest \mac is the one defined for Barnes-Hut treecodes, where $L_j$ is the size of the source cell and $R$ is the distance from the target particle to the center of mass of the source cell; $\theta$, called the opening angle, is the parameter that is used to flexibly control the definition of the well-separated property. In the \fmm framework, a \mac is defined between two cells (rather than a cell and a particle), where $L_i$ is the size of the target cell and $R$ is the distance between the center of mass of the target cell and the center of mass of the source cell.
A flexible definition of the well-separated property cannot be implemented easily in the traditional \fmm framework, where the cell-cell interaction lists are constructed explicitly. This is because one must first provide a list of candidates to apply the \mac to, and increasing the neighbor search domain to $5^3=125$ neighbors instead of $3^3=27$ neighbors and applying the \mac on them is not an efficient solution.
The dual tree traversal is a natural solution to this problem, since the list of candidates for the cell-cell interaction is inherited directly from the parent. This ``inheritance of cell-cell interaction candidates" is lacking from conventional \fmm, and can be provided by the dual tree traversal while simplifying the implementation at the same time. Furthermore, since the dual tree traversal naturally allows interaction between cells at different levels of the tree, it automatically finds the pair of cells that appear in $U,V,W,X$-lists for adaptive trees \cite{LashukETal2009}, but with much greater flexibility (coming from the \mac) and no overhead of generating the lists (since it doesn't generate any).

Another advantage of the dual tree traversal is that it enables the use of non-cubic cells. As we have explained in section \ref{sse:orb}, this  allows using \orb partitioning in the \fmm framework, which has superior load balancing properties compared to Morton ordering and removes the dependence on global indexing. As mentioned above, global Morton indices are a problem when using millions of cores and the tree depth exceeds 21 levels causing the Morton index to overflow from the 64-bit integer. The combination of the dual tree traversal and \orb is an elegant solution to this problem as well.

We describe the dual tree traversal in Algorithm \ref{al:evaluate}, which calls an internal routine for the interaction of a pair of cells, given in Algorithm \ref{al:interact}. First, a pair of cells is pushed to a stack. It is most convenient to start from the root cell although there is no possibility that two root cells will interact. For every step in the while-loop, a pair of cells is popped from the stack and the larger cell is subdivided. Then, Algorithm \ref{al:interact} is called to perform either particle-particle (\PP) or multipole-local (\ML) interactions. If the cells in the pair are too close and either of the cells has children, the pair is pushed to the stack and will be handled later. A more detailed explanation  is given in a previous publication focusing on hybrid treecode \& \fmm techniques \cite{YokotaBarba2012a}.

\begin{algorithm}[t]
\caption{Evaluate()}
\label{al:evaluate}
\begin{algorithmic}
\STATE A=B=rootcell
\STATE push pair(A,B) into a stack
\WHILE{the stack is not empty}
\STATE Pop stack to get a pair(A,B)
\IF{radius of A $>$ radius of B}
\FOR{all child cells ``a" of cell A}
\STATE Interact(a,B)
\ENDFOR
\ELSE
\FOR{all child cells ``b" of cell B}
\STATE Interact(A,b)
\ENDFOR
\ENDIF
\ENDWHILE
\end{algorithmic}
\end{algorithm}

\begin{algorithm}[t]
\caption{Interact(cell A, cell B)}
\label{al:interact}
\begin{algorithmic}
\IF{there are very few particles in both cells, or they don't have child cells}
\IF{particles weren't sent from remote processes}
\STATE Evaluate multipole-local (M2L)
\ELSE
\STATE Evaluate particle-particle (P2P)
\ENDIF
\ELSIF{cells A and B are well separated}
\STATE Evaluate multipole-local (M2L)
\ELSE
\STATE Push pair(A,B) into a stack
\ENDIF
\end{algorithmic}
\end{algorithm}

\subsection{Local essential tree}\label{sse:let}

The partitioning techniques discussed in section \ref{sse:orb} will assign a local portion of the global tree to each process. However, the \fmm evaluation requires information from all parts of the global tree, and the multipole expansions on remote processes must be communicated. Unlike standard domain-decomposition methods, where only a thin halo needs to be communicated, the \fmm requires a global halo. Fortunately, this global halo becomes exponentially coarser as the distance from the target cell increases, as shown in Figure \ref{fig:periodic}. Therefore, the data required to be communicated from far partitions is very small (but never zero). This results in a non-homogeneous \texttt{alltoallv}-type communication of subsets of the global tree.

Once all the subsets of the global tree are communicated between the processes, one can locally reconstruct a tree structure that contains all the information required for the evaluation in the local process. This reconstructed tree is called the \emph{local essential tree} (\LET), and it is a (significantly smaller) subset of the global tree that contains all the information that is necessary to perform the evaluation. Salmon and Warren \cite{WarrenSalmon1992} introduced two key techniques in this area: one for finding which cell data to send, and the other for communicating the data efficiently.

The determination of which data to send is tricky, since each process cannot see what the adaptive tree structure looks like on remote processes. The solution proposed by Salmon and Warren \cite{WarrenSalmon1992} is to use a conservative estimate and communicate a larger portion of the tree than is exactly required. This conservative estimate is obtained by means of a special \mac described in Figure \ref{fig:mac} as the \LET \mac. There, the distance $R$ is defined as the distance between the center of mass of the source cell and the edge of the partition on the remote process. A formula to calculate $R$ can be given as the following expression, where we assume element-wise Boolean operations over array elements (giving zero for false, and 1 for true) and element-wise multiplication and summation:
\begin{eqnarray*}
R&=&|\Delta\mathbf{x}|\\
\Delta\mathbf{x}&=&(\mathbf{x} > \mathbf{x}_{max})(\mathbf{x} - \mathbf{x}_{max}) + 
(\mathbf{x} < \mathbf{x}_{min})(\mathbf{x} - \mathbf{x}_{min})
\end{eqnarray*}

\noindent Here, $\mathbf{x}_{min}$ and $\mathbf{x}_{max}$ are the minimum and maximum coordinate values for all particles in the target partition, while $\mathbf{x}$ is the center of mass of the source cell. This definition of the \LET \mac corresponds to assuming the case where the target cell is located at the edge of the remote partition. We also assume that the target cell is of the same size as the source cell ($L_i = L_j$), so $L_i+L_j$ becomes $2L_j$. These assumptions generally hold quite well and the required part of the tree is sent to the remote process most of the time. However, we must ensure that the \fmm code still works for extreme cases where necessary information fails to be sent. With this end in view, we have added a conditional statement in the interaction calculation (Algorithm \ref{al:interact}) as a further precaution for anomalous cases. This way, the traversal will perform the \ML translation with the smallest cell that is available. This cell may be too large to satisfy the \fmm \mac so there will be a small penalty on the accuracy, but since the occurrence of such a case is so rare, it does not affect the overall result.

\begin{figure*}
\centering
\includegraphics[width=1.0\linewidth]{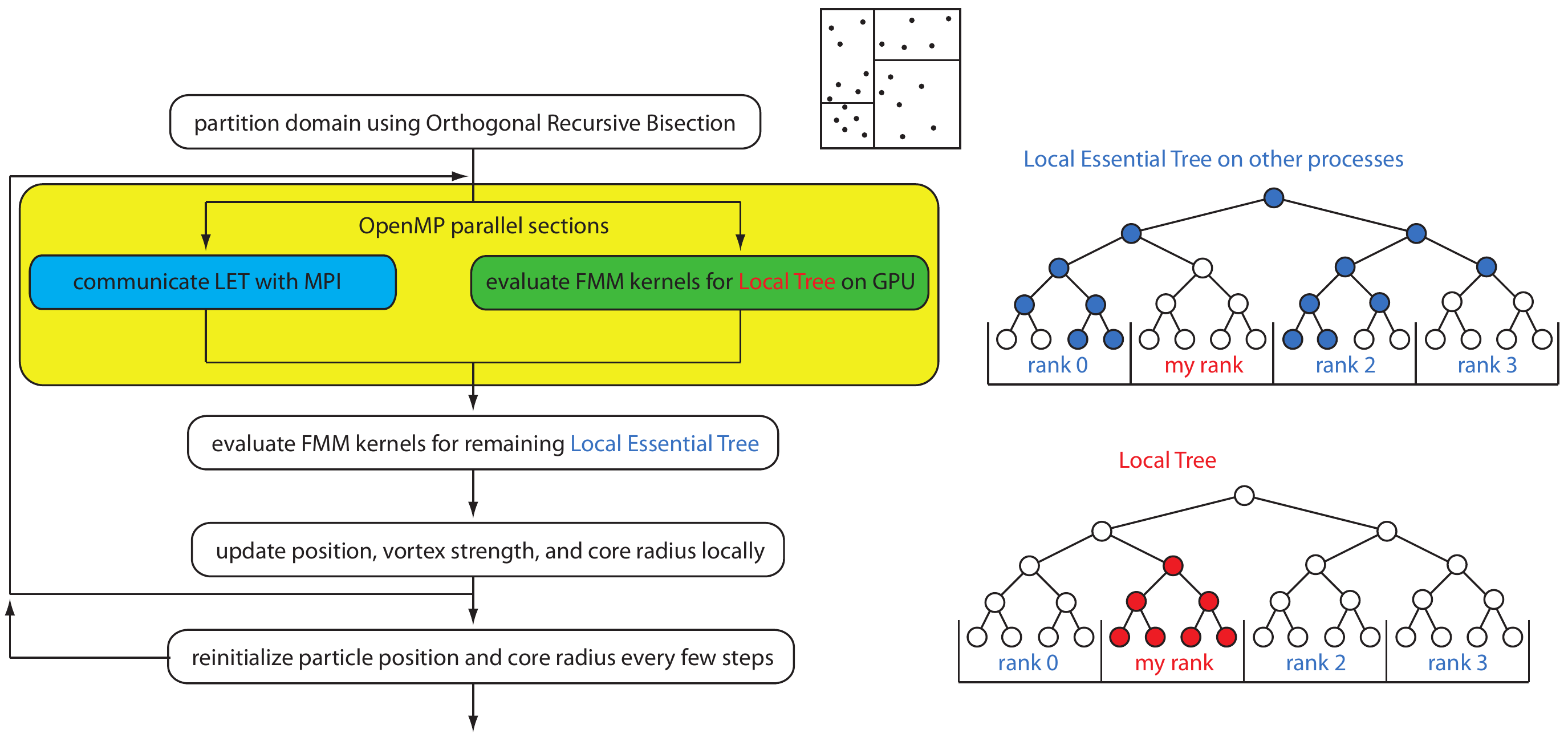}
\caption{Flowchart explaining the use of partitioning and communication strategies of the \fmm in a vortex method application code (left) and schematic of the local essential tree (right).}
\label{fig:flow_chart}
\end{figure*}

A schematic showing how the partitioning and communication techniques are used in the vortex method is given in Figure \ref{fig:flow_chart}. First, the domain is partitioned using the recursive multisection described in section \ref{sse:orb}. Then the \fmm kernels for the local tree are evaluated while the \LET data is being communicated in a separate OpenMP section. After the \LET data is communicated, the \fmm kernels are evaluated again for the remaining parts of the \LET. Subsequently, the position, vortex strength $\boldsymbol{\gamma}$, and core radius $\sigma$ of the vortex particles are updated locally. This information is communicated in the next stage when the \LET is exchanged. In addition, the Lagrangian vortex method needs to reinitialize the particle positions to maintain sufficient overlap of the smooth particles. For this procedure, we can reuse the same tree structure since the particles are reinitialized to the same position every time. Therefore, the partitioning is performed only once at the beginning of this simulation.

\subsection{Periodic FMM}

\begin{figure}[t]
\centering
\includegraphics[width=1.0\linewidth]{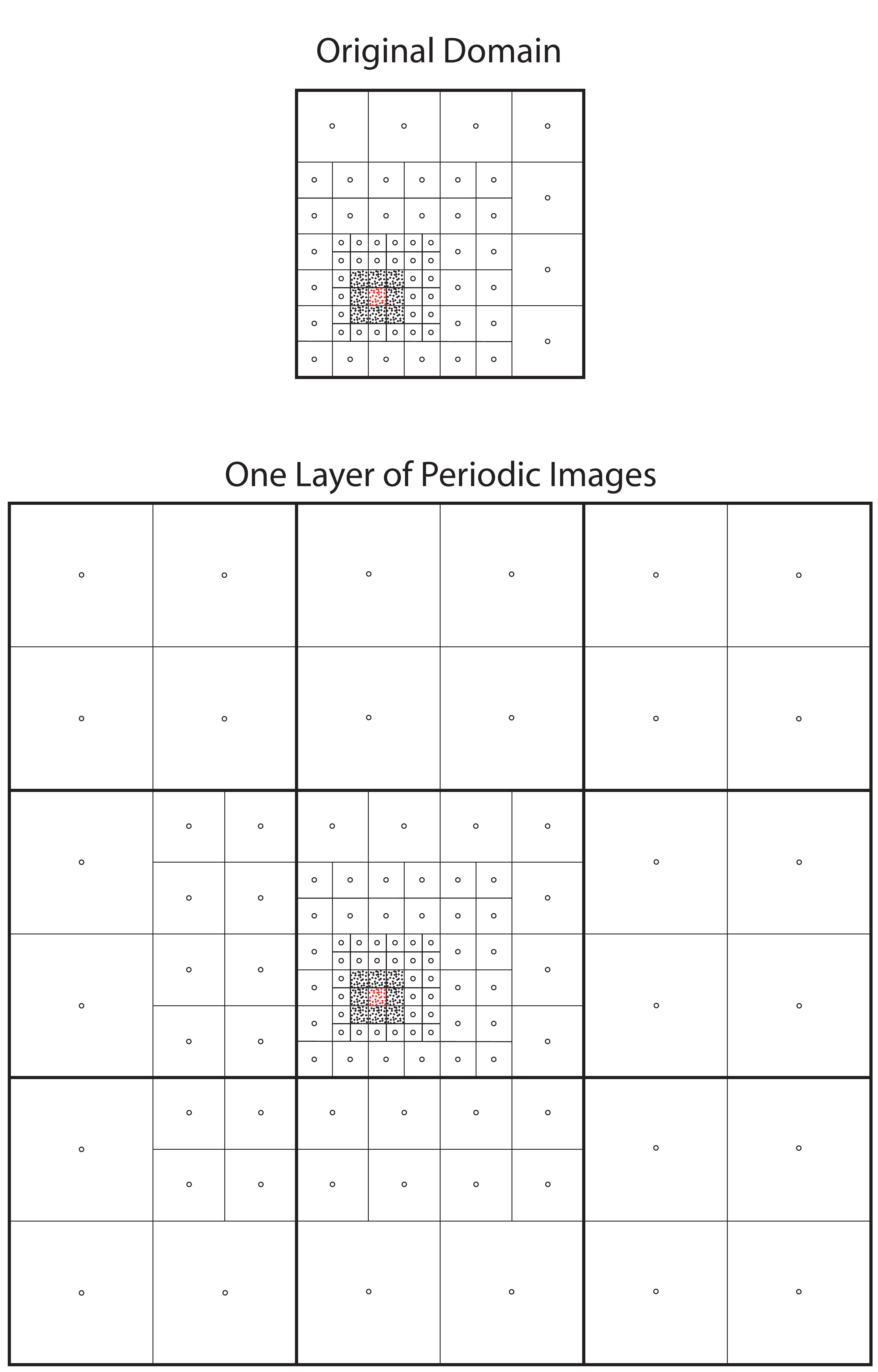}
\caption{Using periodic image domains with the \fmm interaction list, for a target cell illustrated by the red dots.}
\label{fig:periodic}
\end{figure}

The \fmm was originally devised to solve potential problems with free-field boundary conditions. The method can be extended to handle periodic boundary conditions by placing periodic images around the original domain and using multipole expansions to approximate their influence \cite{LambertETal1996}, as illustrated in Figure \ref{fig:periodic}. When a sufficient number of periodic images are placed, the error caused by using a finite number of periodic images becomes smaller than the approximation error of the \fmm itself \cite{LambertETal1996}. This approach to extend the \fmm to periodic domains adds negligible computational overhead to the original \fmm, for two reasons. First, distant domains are clustered into larger and larger cells, so the extra cost of considering a layer of periodic images is constant, while the number of images accounted for grows exponentially. The second reason is that only the sources need to be duplicated and the target points exist only in the original domain. Since the work load for considering the periodicity is independent of the number of particles, it becomes negligible as the problem size increases. An earlier study showed that periodic boundary conditions add approximately 3\% to the calculation time for a million particles \cite{YokotaETal2007}.

\subsection{GPU kernels} \label{sse:gpu}

The \fmm consists of six different computational kernels, as illustrated on Figure \ref{fig:kernels}. In the \exafmm code, all of these kernels are evaluated on \gpu\ devices using \textsc{cuda}. Out of the six kernels, a great majority of the runtime is spent executing \PP and \ML. We use a batch evaluation for these two kernels, since there is no data dependency between them: they are evaluated in one batch after the tree traversal is completed. This batch evaluation can be broken into multiple calls to the \gpu\ device, depending on its storage capacity and the data size. With this approach, we are able to handle problem sizes of up to 100 million particles on a single \gpu, if the memory on the host machine is large enough \cite{YokotaETal2011a}.

\begin{figure}
\centering
\includegraphics[width=1.0\linewidth]{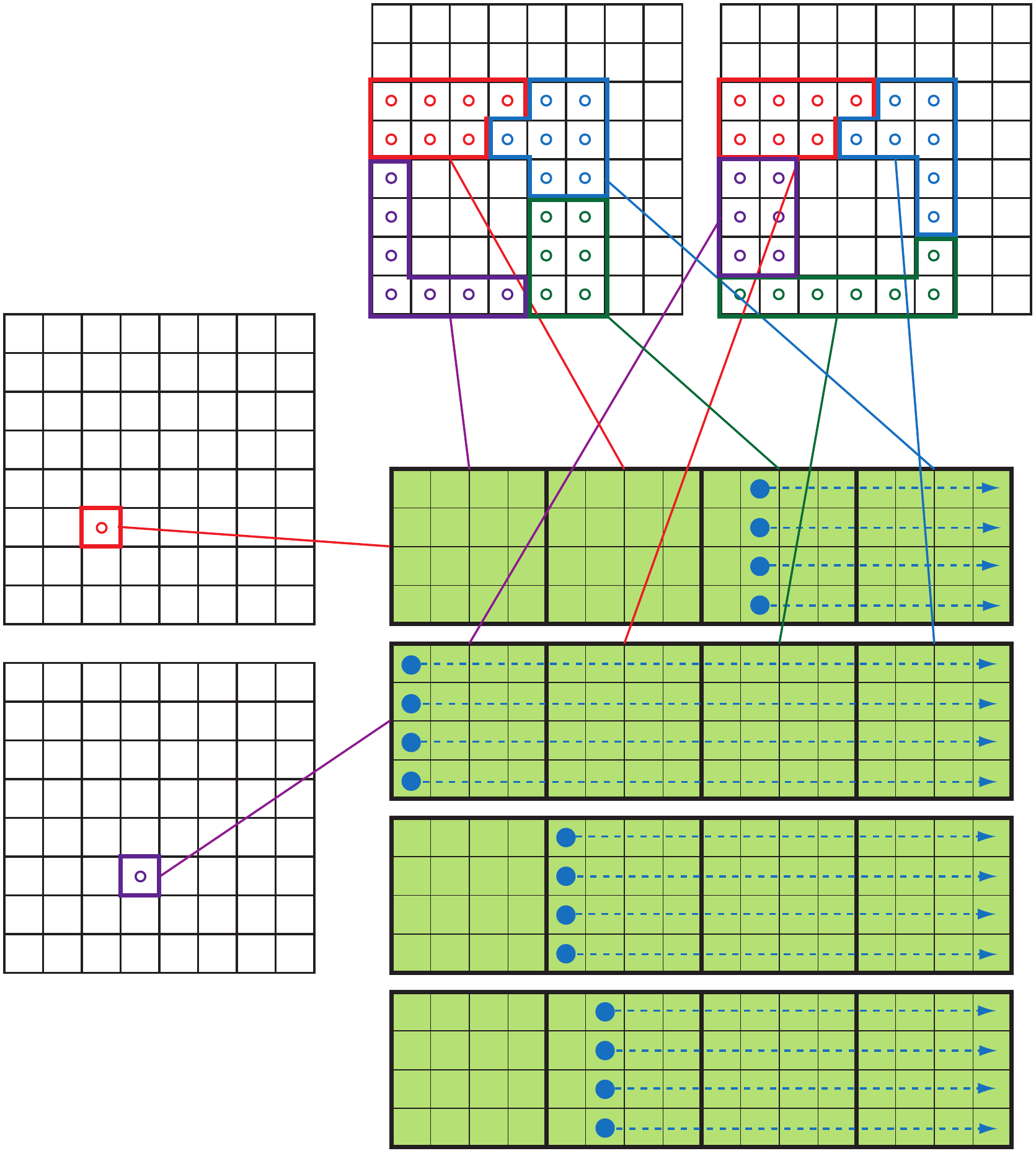}
\caption{Execution of \ML kernel on \gpu s, where each expansion coefficient is mapped to a thread and each target cell is mapped to a thread block. Figure appeared in the chapter by Yokota and Barba \cite{YokotaBarba2010}---permission pending.}
\label{fig:m2l_gpu}
\end{figure}

As an example of the usage of thread blocks in the \gpu\ execution model, we show in Figure \ref{fig:m2l_gpu} an illustration of the \ML kernel on \gpu s.  Each coefficient of the multipole/local expansion is mapped to a thread on the \gpu\ and each target cell is mapped to a thread block, while each source cell is loaded to shared memory and evaluated sequentially. All other kernels are mapped to the threads and thread blocks in a similar manner. More details regarding the \gpu\ implementation of \fmm kernels can be found in chapter 9 of the \emph{GPU Gems Emerald Edition} book \cite{YokotaBarba2010}, with accompanying open-source codes\footnote{In Google Code at \url{http://code.google.com/p/gemsfmm/}}.

\section{Results: parallel simulations of isotropic turbulence}

In this section, we present results from large-scale simulations of homogeneous isotropic turbulence using the \fmm-based vortex particle method, up to a problem size of $4096^3$ computational points. This is still the largest mesh-size for which \dns of turbulence has been published, even though this scale of simulation was achieved 10 years ago. As we discussed in the Introduction, one of the reasons for this is the difficulty in scaling the \fft algorithm beyond a few thousand processes. The current simulations with a vortex method were checked for correctness by comparing the kinetic energy spectrum with that obtained using a trusted spectral-method code (described in \S\ref{s:spectral}). The focus here is not on the physics, however, but on demonstrating large-scale simulation of turbulence using the \fmm as a numerical engine and reporting on performance aspects using \gpu\ hardware. The performance is described via the results of weak scaling tests using between 1 and 4096 processes with $256^3$ (16.8 million) particles per process. Each process offloads to one \gpu\ to speed-up the \fmm, making it even more challenging to obtain good parallel efficiency (it is obviously harder to scale faster code). Despite this, a parallel efficiency of 74\% was obtained for the \fmm-based vortex method on weak scaling between 1 and 4096 processes, with the full application code.  The largest calculation used $69$ billion particles, which as far as we know is the largest vortex-method calculation to date. Previous noteworthy results by other authors  reported calculations with up to 6 billion particles \cite{ChatelainETal2008}, which we surpass here by an order of magnitude.

\subsection{Hardware}

The calculations reported here were run on the \tsubame system during Spring and Fall of 2011, thanks to guest access provided by the Grand Challenge Program of \tsubame. This system  has 1408 nodes, each equipped with two six-core Intel Xeon~{\small X5670}  (formerly Westmere-EP) 2.93GHz processors, three \nvidia~{\small M2050} \gpu s, 54 GB of \textsc{ram}  (96 GB on 41 nodes), and 120 GB of local \textsc{ssd} storage (240 GB on 41 nodes). Computing nodes are interconnected with the InfiniBand device Grid Director 4700 developed by Voltaire Inc., with non-blocking and full bisectional bandwidth. Each node has $2 \times 40$ Gbps bandwidth, and the bisection bandwidth of the system is over 200 Tbps. The total number of~{\small M2050} \gpu s in the system is 4224, and the peak performance of the entire system is 2.4 petaflop/s.

\subsection{Test conditions}

We set up simulations of decaying, homogeneous isotropic turbulence using an initial Taylor-scale Reynolds number of $Re_{\lambda}\approx 500$. Given that we had guest access to the full \tsubame system for a very brief period of time (only a few hours) to produce the scalability results, we opted for a lower Reynolds number than previous turbulence simulations of this size using spectral methods. There is limited experience using vortex methods for \dns of turbulence, but previous work has suggested that higher resolutions are needed than when using the spectral method at the same Reynolds number. We thus decided to be conservative and ensure that we obtained usable results from these one-off runs.

The calculation domain is a box of size $[-\pi,\pi]^3$ with periodic boundary conditions in all directions, using $3^3$ periodic images in each dimension in the periodic \fmm (see Figure \ref{fig:periodic}). To achieve the best accuracy from the \fmm for the present application, the order of multipole expansions was set to $p=14$; this may be a conservative value, but given our limited allocation in the \tsubame system, we did not have the luxury of tuning the simulations for this parameter. 

The \fmm kernels are run in single precision on the \gpu\, which may raise some concerns. We are able to achieve double-precision accuracy using single-precision computations in the \fmm kernels by means of two techniques. The multipole-expansion coefficients have exponentially smaller magnitude with increasing order of expansion; therefore, by adding them from higher- to lower-order terms, we can prevent small numbers being added to the large numbers. This preserves the significant digits in the final sum. The second technique consists of  normalizing the expansion coefficients to reduce the dynamic range of the variables, allowing the use of single-precision variables to get a double-precision result. Thus, we are able to achieve sufficient accuracy to reproduce the turbulence statistics (as shown below) and obtain the same results that would be obtained using double precision (given that the \fmm error is larger than 6 significant digits anyway).

\subsection{Isosurface of the second invariant and kinetic energy spectrum of turbulence}

\begin{figure}
\centering
\includegraphics[width=1.0\linewidth]{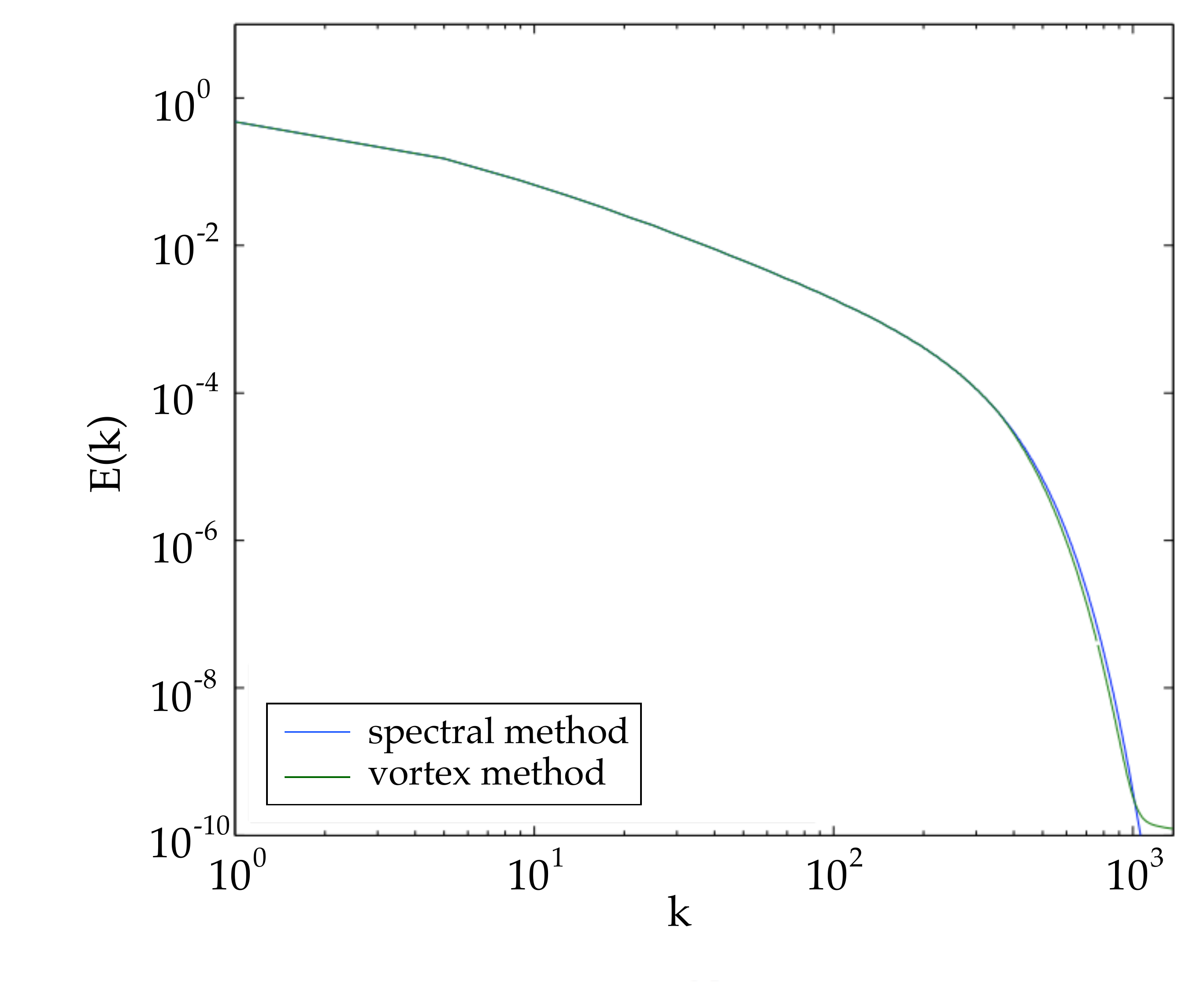}
\caption{Kinetic energy spectrum at $t/T=2$, obtained with the vortex method and spectral method. Notice that the vertical axis goes down to $10^{-10}$, which is many times smaller than in similar plots presented by other authors.}
\label{fig:spectrum}
\end{figure}

\begin{figure*}
\centering
\includegraphics[width=1.0\textwidth]{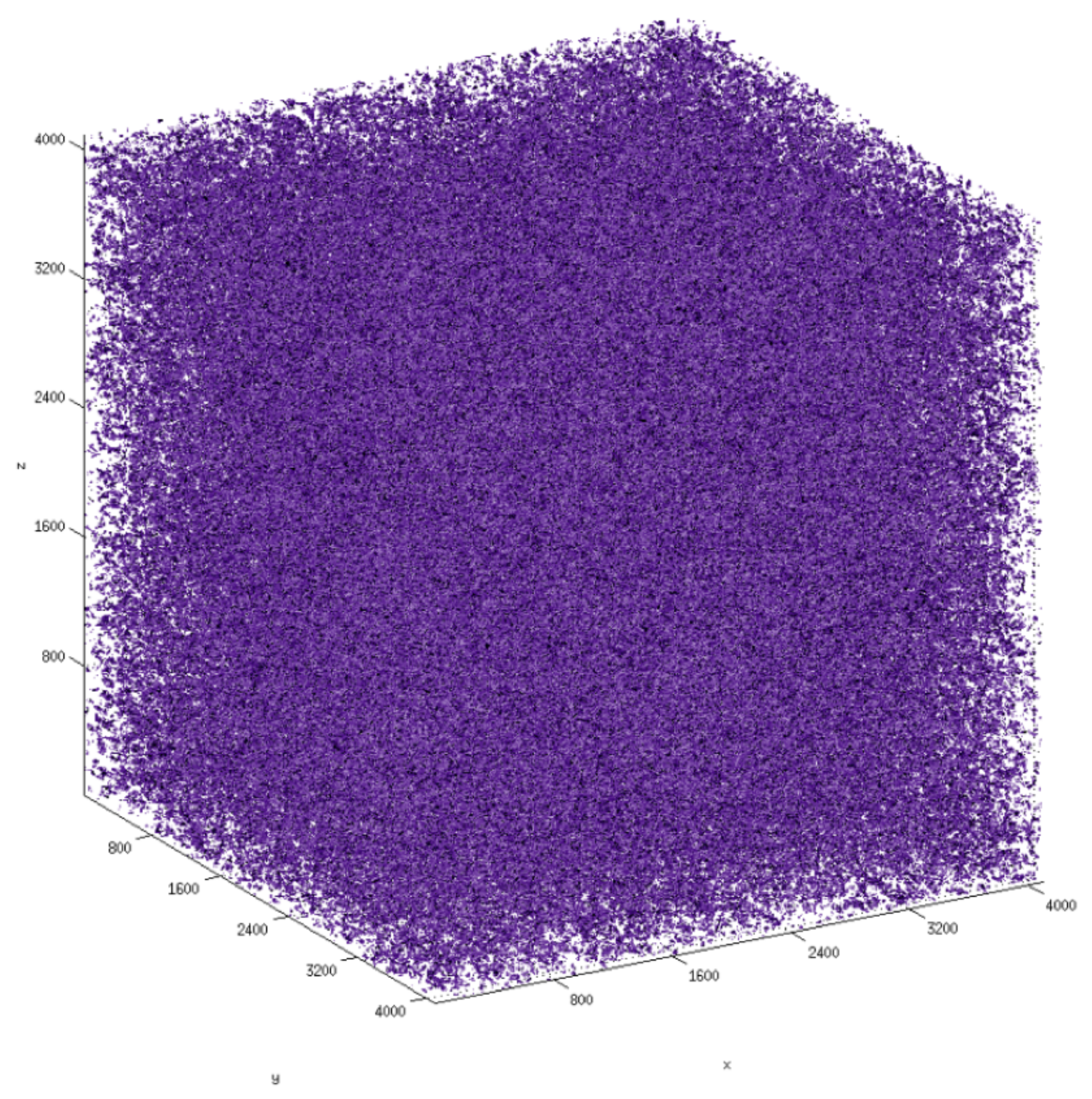}
\caption{Plot of an isosurface of the second invariant of the velocity gradient tensor $II=-10$ for the $4096^3$ mesh at $t/T=2$.}
\label{fig:isosurface}
\end{figure*}

The isosurface of the second invariant of the velocity gradient tensor is shown in Figure \ref{fig:isosurface}. This is a snapshot at the early stages of the simulation and we do not observe any large coherent structures. In order to take a closer look at the quantitative aspects of the vortex simulation,  in Figure \ref{fig:spectrum} we compare the kinetic energy spectrum with that of the spectral method, where $T$ is the eddy turnover time. The energy spectrum of the vortex method is obtained by calculating the velocity field on a uniform lattice, which is induced by the Lagrangian vortex particles. The capability of \texttt{exaFMM} to calculate for different sources and targets enabled such operations. We have excellent quantitative agreement between the vortex method and spectral method and we conclude that our \fmm-based particle method is capable of simulating turbulence of this scale correctly. We turn our attention to the  performance of these calculation in the next section.

\subsection{Weak scaling}

\begin{figure}
\centering
\includegraphics[width=1.0\linewidth]{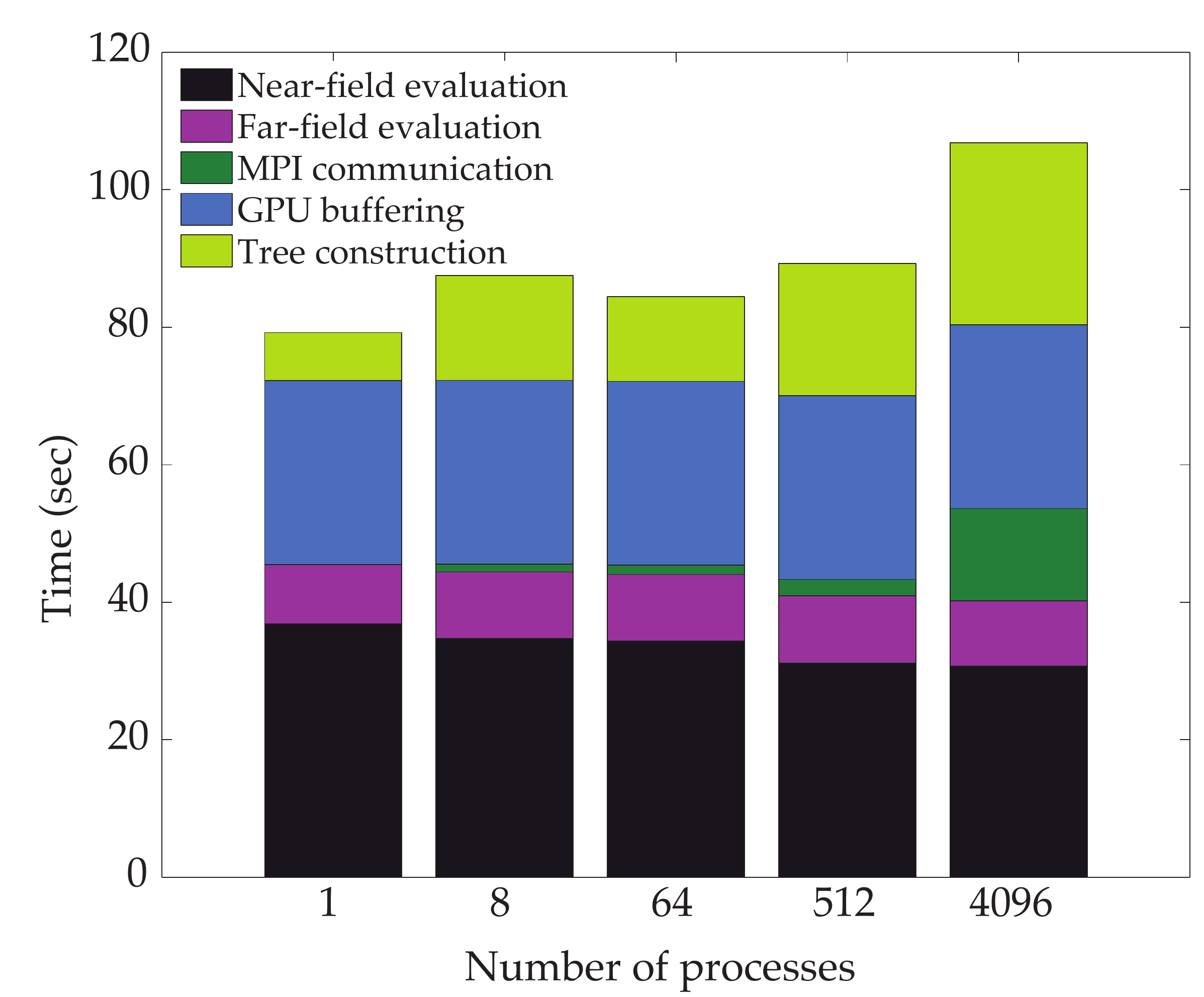}
\caption{Runtime of weak scaling tests, using one \gpu\ per process and three \gpu s per node. All the bars would have the same height if the scaling was perfect. The bars also show a timing breakdown for important phases of the calculation; communication time is overlapped, and is almost completely hidden up to 512 processes.}
\label{fig:weak_scaling}
\end{figure}

We had two very short windows of time in which we were able to run weak scaling tests, one with half the system and the other with almost the full \tsubame system. The larger scaling tests used $256^3$ particles per process, on 1 to 4096 processes and executing three \mpi processes per node, each process assigned to a \gpu\ card within the node.
The results of the weak scaling test are shown in Figure \ref{fig:weak_scaling} in the form of total runtime of the \fmm, and timing breakdown for different phases in the computation. The label `Near-field evaluation' corresponds to the \PP kernel evaluation illustrated in Figure \ref{fig:kernels}, and the label `Far-field evaluation' corresponds to the sum of all the other kernel evaluations, i.e., the far field. The `\mpi communication' is overlapped with the \fmm evaluation (see Figure \ref{fig:flow_chart}), so the plot shows only the amount of time that communication exceeds the local portion of the `near-field' and `far-field' evaluation. In this way, the total height of each bar correctly represents the total wall-clock time of each calculation. Note that particle updates in the vortex-method calculation take less than 0.01\% in all runs and thus were invisible in the bar plots, so we've left this computing stage out of the labels.

The `\gpu\ buffering' label corresponds to the time it takes to form a queue of tasks and corresponding data buffer to be transferred to the \gpu, which is a significant amount of time. We have found this buffering to be necessary in order to achieve high efficiency in the \PP and \ML evaluation and \fmm evaluation on \gpu s. Moreover, this part of the computation scales perfectly, and does not affect the scalability of the \fmm. The parts that do affect the scalability are the tree construction and \mpi communication. Actually, the tree construction also involves \mpi communications for the partitioning, so the parallel efficiency in weak scaling is fully determined by \mpi communications.  Figure \ref{fig:weak_scaling} shows that the current \fmm is able to almost completely hide the communication time up to 512 \gpu s.

It may be worth noting that the $N$-D-hypercube-type communication of the \LET \cite{LashukETal2009} turned out to be slower on \tsubame than a simple call to \texttt{MPI\_Alltoallv} for sending the entire \LET at once.  This is a consequence of the network topology of \tsubame (with a dual-QDR InfiniBand link to each node and non-blocking full-bisection fat-tree interconnect) and also of the relatively small number of \mpi processes. The results shown in Figure \ref{fig:weak_scaling} are those obtained with \texttt{MPI\_Alltoallv} and not the hypercube-type communication.

\begin{figure}
\centering
\includegraphics[width=1.0\linewidth]{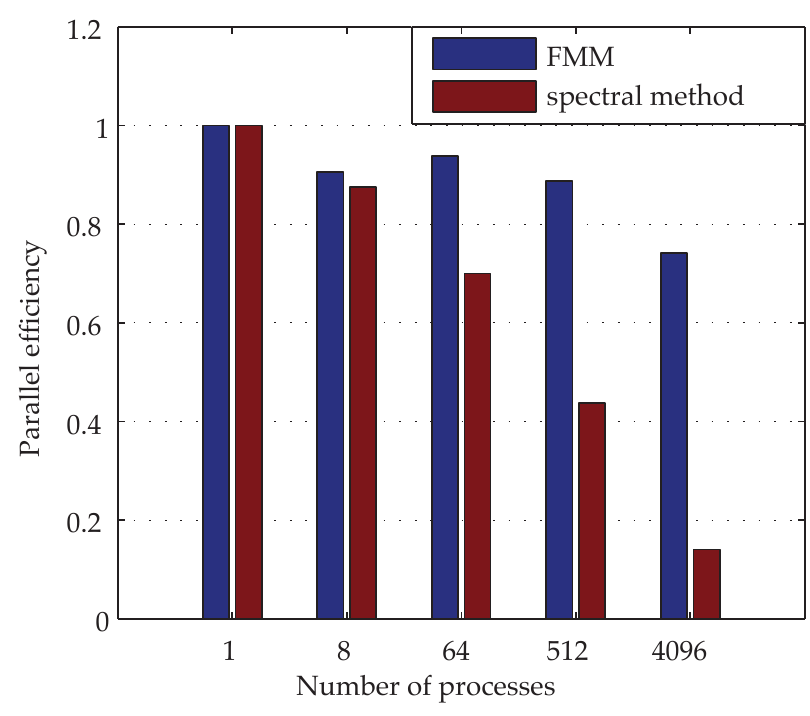}
\caption{Weak scaling from 1 to 4096 processes of the parallel \fmm-based fluid solver on \gpu s, and of the pseudo-spectral method on \cpu s. The parallel efficiency of the \fmm-based solver is 74\% at 4096 processes (one \gpu\ per \mpi process, 3 \gpu s per node). The tests were ran on the \tsubame system in October 2011 using revision 191 of the \texttt{exaFMM} code.  Figure, plotting script and dataset available online and usage licensed under \ccby\cite{YokotaBarba-share92425}.}
\label{fig:compare_scaling}
\end{figure}

We performed a corresponding weak scalability test for the spectral method, increasing the problem size from $256^3$ on one process to $4096^3$ on 4096 processes. We used three \mpi processes per node to match the condition of the vortex method runs, but there is no \gpu\ acceleration in this case. Matching the number of \mpi processes per node should give both methods an equal advantage/handicap for the bandwidth per process. Note that using \gpu s for the \fft within the spectral-method code is unlikely to provide any benefits, because performance improvements of  \texttt{cufft} over \texttt{fftw} would be canceled out by data transfer between host and device\footnote{See \url{http://www.sharcnet.ca/?merz/CUDA_ benchFFT/}} and inter-node communications in parallel runs. Figure \ref{fig:compare_scaling} shows the parallel efficiency obtained with the two methods, under these conditions.

The parallel efficiency of the \fmm-based vortex method is 74\% when going from one to 4096 \gpu s, while the parallel efficiency of the spectral method is 14\% when going from one to 4096 \cpu s. The bottleneck of the spectral method is the all-to-all communication needed for transposing the slabs into pencils as shown in Figure \ref{fig:spectral_fft}. Even though this may not be the best implementation of a parallel \fft, the difference in the scalability between the spectral method and \fmm-based vortex method is considerable. The actual calculation time is in the same order of magnitude for both methods at 4096 processes: it was 108 seconds per time step for the vortex method and 154 seconds per time step for the spectral method. Therefore, the superior scalability of the \fmm has merely closed the gap with \fft, being barely $1.44\times$ faster at this scale. However, we anticipate that this trend will affect the choice of algorithm in the future, as we move to architectures with higher degree of parallelism.

\begin{figure}
\centering
\includegraphics[width=1.0\linewidth]{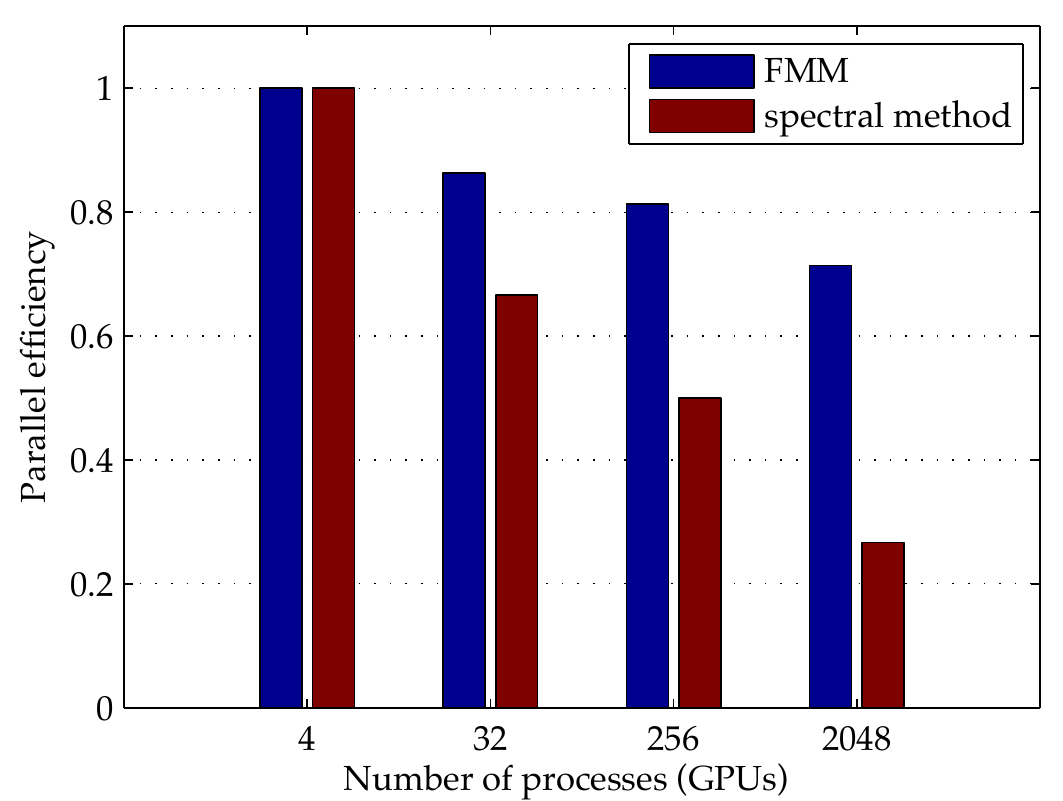}
\caption{Weak scaling from 4 to 2048 processes of the parallel \fmm-based fluid solver on \gpu s, and of the pseudo-spectral method on \cpu s, using 4 million particles per process. The parallel efficiency on 2048 processes is 72 \% for the \fmm-based solver, and 27 \% for the \fft-based solver. These tests were run on the \tsubame system in April 2011, using revision 146 of the \texttt{exaFMM} code.}
\label{fig:compare_scaling_old}
\end{figure}

The scaling test with half the \tsubame system was done several months before and with a different revision of the code, with many changes having been incorporated since then. We include the results here for completeness; see Figure \ref{fig:compare_scaling_old}. In this case, the number of particles per process is much smaller, at 4 million (compared to 16.8 million particles per process in the larger test) and we scale from 4 to 2048 processes. The parallel efficiency of the fmm-based vortex method is 72\% when going from 4 to 2048 \gpu s, while the parallel efficiency of the spectral method was 14\%. When we compare the parallel efficiency with the 1 to 4096 \gpu case, we see that the 4 to 2048 case is scaling relatively poorly. This is due to the number of particles per process being roughly 1/4 in the 4 to 2048 case, and also the improvement in the \texttt{exaFMM} code during the half year gap between the two runs. The run time per time step in this case was 27 seconds for the vortex method and 20 seconds for the spectral method to compute a domain with $2048^3$ points. Note that if we read from the plot in Figure 2 of Donzis et al.~\cite{DonzisETal2008}, their spectral-method calculation on a $2048^3$ grid using 2048 cores of the Ranger supercomputer takes them about 20 seconds per time step. Since this is the same timing we obtained with the \texttt{hit3d} code, we are satisfied that this code provides a credible reference point for the scalability of spectral DNS codes.

\subsection{Sustained performance}

The calculation in the \fmm is mostly dominated by the floating point operations in the particle-particle interactions, while all other parts are a minor contribution in terms of flop/s (although not negligible in terms of runtime). We will thus consider in the estimate of sustained performance only the operations executed by the \PP kernels. Two separate equations are being calculated for the particle-particle interactions: the Biot-Savart equation (\ref{eq:biotsavart}) and the stretching equation (\ref{eq:stretching}). The number of floating point operations required by these two kernels is summarized in Table \ref{tab:flops}.

\begin{table}[h]
\caption{Floating point operations per \PP interaction.}
\label{tab:flops}
\begin{center}
\begin{tabular}{|c|c|c|}
\hline
Operation & Biot-Savart & Stretching \\
\hline
+ & 19 & 25 \\
\hline 
- & 14 & 18 \\
\hline
* & 32 & 56 \\
\hline
/ & 1 & 1 \\
\hline
sqrtf & 1 & 1 \\
\hline
rsqrtf & 1 & 1 \\
\hline
expf & 2 & 2 \\
\hline
Total & 70 & 104 \\
\hline
\end{tabular}
\end{center}
\end{table}

The approximate number of flop/s for one step of the vortex method calculation of isotropic turbulence is obtained by the following equation.
\begin{eqnarray*}
flop/s &=&(\#\ processes) \\
&\times&(target\ particles\ per\ process) \\
&\times&(source\ cells\ per\ target) \\
&\times&(source\ particle\ per\ cell) \\
&\times&(flop/s\ per\ interaction) \\
& / &(wall\ clock\ time) \\
&=&4096\times(6.8\times10^7)\times19\times512\times174/108 \\
&=&1.08\times10^{15}=\mathbf{1.08\ petaflop/s}
\end{eqnarray*}

\noindent Thus, the current \fmm-based vortex method achieved a sustained performance of 1.08 petaflop/s (single precision) on 4096 \gpu s of \tsubame.

\subsection{Reproducibility and open-source policy}

The authors of the \exafmm code have a consistent policy of making science codes available openly, in the interest of reproducibility. The entire code that was used to obtain the present results is available from \url{https://bitbucket.org/exafmm/exafmm}. The revision number used for the results presented in this paper is 191 for the large-scale tests up to 4096 \gpu s. Documentation and links to other publications are found in the project homepage at \url{http://exafmm.org/}. Figure \ref{fig:compare_scaling}, its plotting script and datasets are available online and usage is licensed under \ccby\cite{YokotaBarba-share92425}.

We acknowledge the use of the \texttt{hit3D} pseudo-spectral DNS code for isotropic turbulence, and appreciate greatly their authors for their open-source policy; the code is available via Google code at \url{http://code.google.com/p/hit3d/}.

\section{Conclusions}

This work represents several milestones. Although the \fmm algorithm has been taken to petascale before (notably, with the 2010 Gordon Bell prize winner), the present work represents the first time that this is done on \gpu\ architecture. Also, to our knowledge, the present work is the largest direct numerical simulation with vortex methods to date, with 69 billion particles used in the cubic volume; this is an order of magnitude larger than the previously reported record. Yet another significant event is reaching a range where the highly scalable \fmm starts showing advantage over \fft-based algorithms. With a 1.08 petaflop/s (single precision) calculation of isotropic turbulence in a $4096^3$ box, using 4096 \gpu s, we are within reach of a turning point. The combination of application, algorithm, and hardware used are also notable.

The real challenge in exascale computing will be the optimization of data movement. When we compare the data movement of \fmm against other fast algorithms like multigrid and \fft, we see that the \fmm has a potential advantage. Compared to \fft, both multigrid and \fmm have an advantage in the asymptotic complexity of the global communication. The hierarchical nature of multigrid and \fmm results in $\mathcal{O}(\log P)$ communication complexity where $P$ is the number of processes. On the other hand, \fft requires two global-transpose communications between $\sqrt{P}$ processes, and has communication complexity of $\mathcal{O}(\sqrt{P})$. When $P$ is in the order of millions, it seems unrealistic to expect that an affordable network can compensate for this large gap in the communication complexity. Although it is not the focus of the present article, we would like to briefly note that an advantage of \fmm over multigrid is obtained from differences in the synchronization patterns. For example, increasing the desired accuracy in iterative solvers using multigrid will result in more iterations, hence more global synchronizations. Conversely, increasing the accuracy in \fmm involves increasing the number of multipole expansions, which results in even higher arithmetic intensity in the inner kernels while the number of global synchronizations remains the same. As the amount of concurrency increases, bulk-synchronous execution/communication models are reaching their limit. Thus, \fmm provides a new possibility to reduce the amount of communication and synchronization in these inherently ``global" problems.

Finally, we would like to point out that the \fmm can be used to solve the Poisson equation directly, or as a preconditioner for an iterative solver. Therefore, we are not concerned at this point about the fact that vortex methods may still be comparatively inefficient for the simulation of fluid turbulence, where spectral methods will continue to dominate in the foreseeable future. This does not detract from the conclusions about the efficiency of the \fmm itself, which is the object of our study. In future work, we would like to demonstrate the efficiency of \fmm by using it as a Poisson solver or preconditioner in the framework of more standard finite difference/volume/element methods. There, the comparison against \fft and multigrid methods should be of interest to a broader spectrum of the \cfd community.

\section{Acknowledgments}
Computing time in the \tsubame system was made possible by the Grand Challenge Program of \tsubame. The current work was partially supported by the Core Research for the Evolution Science and Technology (CREST) of the Japan Science and Technology Corporation (JST). LAB acknowledges funding from NSF grant OCI-0946441, ONR grant \#N00014-11-1-0356 and NSF CAREER award OCI-1149784. LAB is also grateful for the support from \nvidia Corp.\ via an Academic Partnership award (Aug.~2011).

\end{document}